\documentclass[a4paper, 10pt]{article}

\usepackage{amsmath}
\usepackage{amsfonts}
\usepackage{hyperref}
\usepackage{authblk}

\usepackage{pgfplots} 				
\usepackage{tikz}
\usepackage{tikzscale}
\usepackage{environ}

\usepackage{multirow}

\graphicspath{{./}}

\newcommand{\R}{\mathbb{R}}
\newcommand{\N}{\mathbb{N}}
\newcommand{\ns}{\mathcal H_k (\Omega)}

\DeclareMathOperator*{\argmax}{arg\max}

\newtheorem{theorem}{Theorem}

\title{Biomechanical surrogate modelling using stabilized vectorial greedy kernel methods}
\author[1]{Bernard Haasdonk \thanks{haasdonk@mathematik.uni-stuttgart.de}}
\author[1]{Tizian Wenzel \thanks{tizian.wenzel@mathematik.uni-stuttgart.de}}
\author[2]{Gabriele Santin \thanks{gsantin@fbk.eu, \href{http://orcid.org/0000-0001-6959-1070}{orcid.org/0000-0001-6959-1070}}}
\author[3]{Syn Schmitt \thanks{schmitt@simtech.uni-stuttgart.de}}

\affil[1]{Institute for Applied Analysis and Numerical Simulation, University of Stuttgart, Germany}
\affil[2]{Center for Information and Communication Technology, Fondazione Bruno Kessler, Italy}
\affil[3]{Institute for Modelling and Simulation of Biomechanical Systems, University of Stuttgart, Germany}

\begin{document}
\newlength\fwidth

\maketitle

\begin{abstract}
Greedy kernel approximation algorithms are successful techniques for sparse and accurate data-based modelling and function approximation. Based on a recent 
idea of stabilization \cite{wenzel2019novel} of such algorithms in the scalar output case, we here consider the vectorial extension built on VKOGA 
\cite{Wirtz2013}. We introduce the so called $\gamma$-restricted VKOGA, comment on analytical properties and present numerical evaluation on data from a 
clinically relevant application, the modelling of the human spine. The experiments show that the new stabilized algorithms result in improved accuracy and 
stability over the non-stabilized algorithms.
\end{abstract}

\section{Introduction} \label{sec:introduction} 

Kernel methods are used in various fields of machine learning or pattern analysis. They yield efficient and flexible ways to recover functions from data since 
they can deal with arbitrarily scattered points. The combination of their flexibility with the strong mathematical theory about e.g.\ existence, convergence, 
stability make them a nice tool for applications \cite{Fasshauer2015, Wendland2005}. \\
In this paper we apply a recently introduced idea that has lead to a new class of stabilized greedy kernel algorithms \cite{wenzel2019novel}, extend it to 
vectorial function approximation and apply it to a real life setting from research in biomechanics. Some theoretic statements can be extended from the scalar to 
the vectorial case. All in all these stabilized methods provide further flexibility and are able to efficiently mitigate the problem of having numerical 
instabilities. \\
The paper is organized as follows. To begin with we recall in Section \ref{sec:stabilized_vkoga} some basics about kernel interpolation with a focus on greedy 
kernel approximation and explain the stabilized extension. Section \ref{sec:appl_spine} gives background information about our application settings and the use 
of kernel methods. The following Section \ref{sec:numerical_exp} explains the conducted numerical experiments as well as the practical results. Section 
\ref{sec:conclusion} concludes with a summary and an outlook.

\section{Stabilized VKOGA algorithm} \label{sec:stabilized_vkoga} 

We start with a nonempty set $\Omega \subset \R^d$. A real-valued kernel is a symmetric function $k: \Omega \times \Omega \rightarrow \R$. For arbitrary points 
$X_N := \{ x_1, .., x_N \} \subset \Omega$ the kernel matrix $A \in \R^{N \times N}$ is a symmetric matrix with entries $A_{ij} = k(x_i, x_j)$. If this kernel 
matrix is positive semi-definite for any set of points $X_N \subset \Omega$, then the kernel is called positive definite. If the kernel matrix is even positive 
definite for any set of pairwise distinct points, then the kernel is called strictly positive definite. In the following we will focus on this case of 
strictly positive definite kernels and we refer to the monographs \cite{Fasshauer2015, Wendland2005} for more details. \\
For any such kernel there is a unique Hilbert space of functions, namely the native space $(\ns, (\cdot, \cdot)_{\ns})$, which is a Reproducing Kernel Hilbert 
Space (RKHS). 
A popular choice is given by radial basis function kernels, i.e.\ the kernel can be expressed with the help of some function $\Phi$ and a kernel parameter 
$\epsilon \in \R$ as $k(x, y) = \Phi(\epsilon \Vert x-y \Vert)$. Examples are given by the Gaussian kernel $\Phi_{\text{Gauss}}(r) = \exp(-(\epsilon \cdot 
r)^2)$and the linear Mat\'ern kernel $\Phi(r) = (1+r) \cdot \exp(-r)$. The decay of the Fourier transform of those radial basis functions is 
decisive for their properties. The Fourier transform of the Gaussian decays exponentially, whereas the Fourier transform of the linear Mat\'ern decays only 
algebraically.

In such RKHS the interpolation of functions - or more general data based approximation tasks - can be analyzed. For a given function $f \in \ns$ and some 
interpolation points $X_N$ the interpolant $s_N$ is given by the orthogonal projection $\Pi_{V(X_N)}(f)$ of $f$ onto $V(X_N) := \{ k(\cdot, x_i), x_i \in X_N 
\}$ and thus can be expressed as
\begin{align*}
s_N(\cdot) = \Pi_{V(X_N)}(f) = \sum_{i=1}^N \alpha_i k(\cdot, x_i), \;\; \alpha_i\in\R, 1\leq i\leq N.
\end{align*}
In some applications the data is affected by noise, so it does not make sense to interpolate the given values, while it is rather advisable to approximate them 
while taking some regularization into account. For this one can consider minimizing $\sum_{i=1}^N \Vert f(x_i) - s_N(x_i) \Vert_2^2 + \lambda \cdot \Vert s_N 
\Vert_{\ns}^2$ which corresponds to solving the linear system
\begin{align}
(A + \lambda \cdot I ) \alpha = y
\label{eq:reg_system}
\end{align}
with $y = (f(x_i))_{i=1}^N$. To measure the interpolation error $\Vert f - \Pi_{V(X_N)}(f) \Vert_{L^{\infty}}$ one can introduce the Power function $P_{X_N}: 
\Omega \rightarrow \mathbb{R}$ as
\begin{align}
\label{eq:power_function_via_sup}
P_N(x) := P_{X_N}(x) = \sup_{0 \neq f \in \ns} \frac{|f(x) - \Pi_{V(X_N)}(f)(x)|}{\Vert f \Vert_{\ns}}.
\end{align}
From this definition we can directly conclude 
\begin{align*}
|f(x) - \Pi_{V(X_N)}(f)(x)| \leq& P_N(x) \cdot \Vert f \Vert_{\ns}.
\end{align*} 

For the analysis of the kernel interpolants geometric quantities about the distribution of the interpolation points are important. The fill distance $h_N$ and 
the separation distance $q_N$ are defined as
\begin{equation}\label{eq:fill_and_sep_dist}
h_N := \sup_{x \in \Omega} \min_{x_i \in X_N} \Vert x-x_i \Vert_2, \hspace{1cm} q_N := \min_{x_i \neq x_j \in X_N} \Vert x_i - x_j \Vert_2.
\end{equation}

A priori it is unclear how to select good interpolation points for a given set of data or some functions. To circumvent this, one applies greedy methods which 
start with an empty set $X_0 = \{ \}$ and iteratively add another interpolation point according to some selection criterion, $X_{N+1} = X_N \cup \{ x_{N+1} 
\}$. 
\\
There are three main selection criteria in the literature, namely $f$-greedy, $f/P$-greedy and $P$-greedy \cite{SchWen2000, DeMarchi2005, Mueller2009} which 
choose the next point from $\Omega$ according to some indicator. For the vectorial case, for $x\in\Omega$ they are:
\begin{enumerate}
\item $f$-greedy: \hspace{7.75mm} $\eta_f^{(N)}(x) = \Vert f(x) - \Pi_{V(X_N)}(f)(x) \Vert_2$
\item $P$-greedy: \hspace{7mm} $\eta_P^{(N)}(x) = P_{X_N}(x)$
\item $f/P$-greedy: \hspace{3.4mm} $\eta_{f/P}^{(N)}(x) = \Vert f(x) - \Pi_{V(X_N)}(f)(x) \Vert_2 /P_{X_N}(x)$.
\end{enumerate}
In order to create a scale of selection criteria which lie in between those known criteria, one introduces a restriction parameter $\gamma \in (0, 1]$ and a 
restricted set $\Omega_{\gamma}^{(N)} := \{ x \in \Omega, P_N(x) \geq \gamma \cdot \Vert P_N \Vert_\infty \}$ and chooses the next interpolation point within 
$\Omega_{\gamma}^{(N)}$ according to some standard selection criterion. This works since the Power function is scalar valued and the interpolation points are 
shared among all dimensions.\footnote{This corresponds to the case of using separable matrix-valued kernels, i.e. $K(x,y) := k(x,y) \cdot I$ where $I$ is the $d 
\times d$ identity matrix \cite{Wittwar2018}. }
For $\gamma = 1$ it holds $\Omega_\gamma^{(N)} = \{ x \in \Omega, P_N(x) = \Vert P_N \Vert_\infty \}$, thus we obtain the standard $P$-greedy algorithm for any 
selection criterion $\eta^{(N)}(x)$. For 
$\gamma = 0$ it holds $\Omega_\gamma^{(N)} = \Omega$, thus we obtain the unrestricted algorithm. The naming \textit{restricted} is obviously related to the 
restriction of the set $\Omega$ to $\Omega_\gamma^{(N)}$, the name \textit{stabilized} is motivated by part of the results within \cite{wenzel2019novel}, which 
are summarized in Theorem \ref{th:summary_stabilized}. If the maximum within a selection rule is not unique, any point realizing the maximum can be picked. 

As an example, the $\gamma$-stabilized $f$-greedy chooses the next point according to
\begin{align*}
x_{N+1} = \argmax_{x \in \Omega_\gamma^{(N)}} \Vert f(x) - \Pi_{V(X_N)}(f)(x) \Vert_2.
\end{align*}

Several rigorous analytical statements can be derived for this kind of algorithms and we will summarize a few of them in the following Theorem 
\ref{th:summary_stabilized}. The proofs are straightforward consequences of those which can be found in \cite{wenzel2019novel}.

\begin{theorem} \label{th:summary_stabilized}
Assume that $\Omega \subset \mathbb{R}^d$ is a compact domain which satisfies an interior cone condition and has a Lipschitz boundary. Suppose that $k$ is a 
translational invariant kernel such that its native space is norm equivalent to the Sobolev space $H^\tau(\Omega)$ with $\tau > d/2$. Then any 
$\gamma$-stabilized algorithm applied to a function in $f\in\ns$ gives a sequence of 
point sets $X_N \subset \Omega$ such that it holds:

\begin{itemize}
\item Lower and upper bound on the Power function:
\begin{align*}
c_P \cdot N^{\frac{1}{2}-\frac{\tau}{d}} \leq \Vert P_{N} \Vert_{L^\infty(\Omega)} \leq C_P \cdot \gamma^{-2} \cdot N^{\frac{1}{2}-\frac{\tau}{d}}.
\end{align*}
\item Asymptotic uniform point distribution:
\begin{align*}
\rho_{X_N} := \frac{h_N}{q_N} \leq c \cdot \gamma^{-4} ~~ \forall N \in \N.
\end{align*}
\item Lower and upper bounds on the smallest eigenvalue:
\begin{align*}
c_1 \cdot \gamma^{8\tau - 4d -4} \cdot N^{1-2\tau/d} \leq \lambda_{\min}(X_N) \leq c_2 \cdot \gamma^{-4} \cdot N^{1-2\tau/d}.
\end{align*}
\end{itemize}
\end{theorem}

Finally we want to point out that the $\gamma$--parameter only affects the choice of points, i.e.\ it modifies the greedy selection. However, if given points 
are used the $\gamma$ parameter does not change the computed interpolant anymore. This is in contrast to the parameters $\lambda$ and $\epsilon$, which modify 
also the interpolant if points are given. \\

From Theorem \ref{th:summary_stabilized} we can conclude that the product $\lambda_{\min} \cdot \Vert P_N \Vert_\infty^2$ is bounded from both sides for $\gamma 
> 0$. This motivates to take the value of the Power function of the previously chosen point as a measure of stability. This will be used in Equation 
\eqref{eq:stopping_by_stability} to implement a stopping criterion based on stability.

\section{Application to spine modelling} \label{sec:appl_spine}

\begin{figure}[ht]
 	\setlength\fwidth{.4\textwidth}
 	\centering
 	\includegraphics[width=0.29\textwidth]{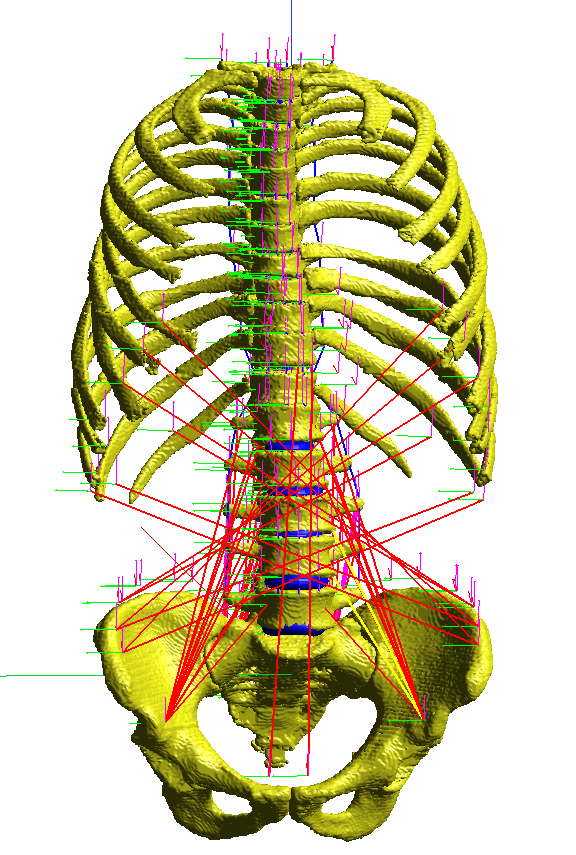}
 	\includegraphics[width=0.65\textwidth]{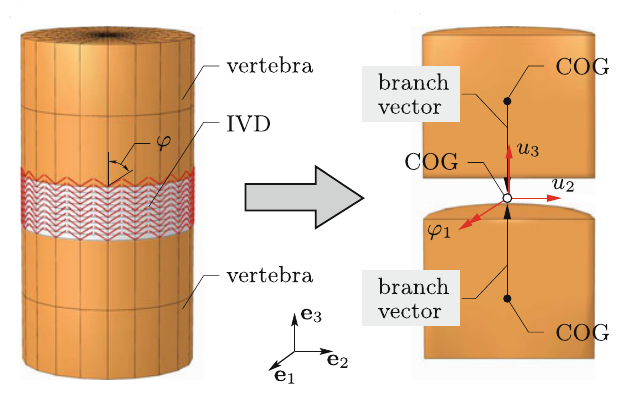}
 	\caption{Visualization of the biomechanical model. On the left the whole spine model is depicted, on the right the modelling scheme of an IVD reduced 
to 
a 3-d force/torque element is shown \cite{karajan2013a}.}
\label{fig:ivd}
\end{figure}

Biomechanical models of the human spine account for the most significant structures which carry the load of daily life. That are mostly the ligaments, the 
muscles with both passive and active contributions and the intervertebral discs (IVDs), of course. An IVD, in this sense, can be seen as a combination of both 
the defining structure for the degrees of freedom between two vertebral bodies and the force and rotational moment transducing elements between these two bones, 
cf. Figure \ref{fig:ivd}.

Mostly, IVDs are modelled by a linear approximation of forces and rotational moments calculated according to the respective displacements, e.g.\ 
\cite{monteiro2011}. Alternatively, as long as quasi-static movements are studied, very detailed, finite element models of isolated IVDs or a combination of few 
spinal segments are used \cite{DREISCHARF20141757}. Another approach to model a reduced IVD used a polynomial approximation and showed that the classical linear 
approximations overestimate actual stiffnesses in the working range \cite{karajan2013a,Rupp2015a}. This observation gave rise to the idea of looking into an 
even more sophisticated mapping of displacements on the input to output forces and rotational moments.
input-output mapping in three spacial dimensions. 

Obviously, kernel modelling seems to be an ideal approach for this need. First, kernel surrogates promise to capture the mapping characteristics well, second, 
extensions to higher input and output dimensions seem feasible and third, compared to respective detailed finite element models surrogate models evaluate the 
mapping stunningly fast \cite{wirtz2015}.

Assuming symmetry in the sagittal and frontal plane, an input-output relation $f: \R^3 \rightarrow \R^3$ is considered and studied, here.

\section{Numerical experiments} \label{sec:numerical_exp}

The considered dataset consists of 1370 input points in $\R^3$ with corresponding output points in $\R^3$. 1238 points are used for training and validation, 
the 
remaining points are used as a test set. No scaling is applied to the data. In order to show the flexibility and thus improved accuracy of the stabilized 
algorithms on the presented data set, we compute unstabilized approximants, used as base models, as well as stabilized models for the $f$- and the $f/P$-greedy 
algorithm. For both the unstabilized and stabilized models we also use regularization in a second step. The experiments are related to those in 
\cite{SH_bookchapter, wirtz2015}, however due to different setups they are not identical.  \\
The base models are given by standard kernel surrogates where the point selection is done either with vectorial $f$-greedy or $f/P$-greedy. To evaluate good 
parameters, first of all a $5$-fold cross validation is run on $20$ logarithmic equally spaced kernel parameters $\epsilon$. The best $\epsilon$ value is used 
for a second step, where the best $\lambda$ parameter from Equation \eqref{eq:reg_system} is evaluated with help of another $5$-fold cross validation. For this 
we use $20$ logarithmic equally spaced values between $10^{-16}$ and $10^{3}$. As an error measure we use the Root Mean Square Error (RMSE)

\begin{align*}
E_{\text{RMSE}}(s, X, Y) = \left( \frac{1}{|X|} \cdot \sum_{i=1}^{|X|} \Vert s(x_i) - y_i \Vert_2^2 \right)^{1/2}.
\end{align*}
The stabilized models are given by kernel surrogates where the points are selected with help of stabilized vectorial $f$- or $f/P$-greedy algorithms. We start 
by using the same kernel parameter $\epsilon$ which was selected for the base model and run instead a $5$-fold cross validation on $11$ equally spaced 
stabilization parameters $\gamma \in [0, 1]$. As a second step we evaluate again the best $\lambda$ parameter with help of a $5$-fold cross validation. This 
procedure keeps the computation time similar to the base model.

The used hyperparameters are summarized in Table \ref{tab:hyperparameters}. For the experiments the linear Mat\'ern kernel is used since it satisfies all the 
prerequisites of Theorem~\ref{th:summary_stabilized}. 

\begin{table}
\centering
\begin{tabular}{|c|ccc|ccc|ccc|} 
\hline
$k$ & $\epsilon_{\min}$	& $\epsilon_{\max}$	& $n_{\epsilon}$ & $\gamma_{\min}$	& $\gamma_{\max}$	& $n_{\gamma}$ & $\lambda_{\min}$	& 
$\lambda_{\max}$	& $n_{\lambda}$ \\ \hline
5 & $10^{-2}$ & $10^1$ 	& $20$ 	& $0$ & $1$ & $11$ & $10^{-16}$ & $10^3$ & 20 \\
\hline
\end{tabular}
\caption{Overview of the hyperparameter ranges. The $\gamma$ values are equally spaced, the others are logarithmically equally spaced.}
\label{tab:hyperparameters}
\end{table}

The greedy selection algorithms stop either when all points within the training set are selected or if some threshold on the residual or on the Power function 
is met. As a tolerance on the residual we use $\tau_f = 10^{-7}$, that means the selection is stopped if $\max \Vert s_N(x_i) - y_i \Vert_2 < \tau_f$ is met. 
As 
a 
tolerance for the Power function we use $\tau_P = 10^{-3}$ and the selection is stopped if 
\begin{align}
P_N(x_{N+1}) < \tau_P.
\label{eq:stopping_by_stability}
\end{align}
We remark that this last criterion is directly linked to the stability. If points with small Power function value are selected, it means that the interpolation 
points cluster. We recall that this means in particular that $\lambda_{\min}(X_{N})$ is below a certain threshold, making further computations unstable. 
Moreover, although a thorough discussion on the fine tuning of these thresholds is beyond the scope of this paper, we remark that the chosen values 
appear to be reasonable in this setting since they are sufficient to achieve the desired accuracy, while avoiding instabilities.

Table \ref{tab:num_results} lists both the hyperparameters which were selected by the cross-validations and the resulting accuracies of the interpolants. The 
$E_{\max, \text{rel}}$ and the $E_{\text{RMSE}, \text{rel}}$ errors are defined according to
\begin{align*}
E_{\max, \text{rel}} &:= \max_{i=1, .., |X|} \Vert s(x_i) - y_i \Vert_2 / \Vert y_i \Vert_2,\\
E_{\text{RMSE}, \text{rel}} &:= \left( \frac{1}{|X|} \cdot \sum_{i=1}^{|X|} \frac{\Vert s(x_i) - y_i \Vert_2^2}{\Vert y_i \Vert_2^2} \right)^{1/2}.
\end{align*}

\begin{table}
\hspace*{-1.2cm}\small{
\centering
\begin{tabular}{|c|l|l|c|l|l|l|} 
\hline
			& \multicolumn{2}{|c|}{$f$-greedy} 					& & \multicolumn{2}{|c|}{$f/P$-greedy} \\ \hline
			& \multicolumn{1}{|c|}{Hyperparameters}				& \multicolumn{1}{|c|}{Results} && 
\multicolumn{1}{|c|}{Hyperparameters} 
& \multicolumn{1}{|c|}{Results}				\\ \hline
\multirow{4}{*}{\rotatebox[origin=c]{90}{\textbf{base}}}		& $\epsilon_{\text{base}} \hspace{.4mm} = 6.158\cdot 10^{-2}$ 	&	$E_{\max} 
\hspace{5.1mm} = 347.22$ 	&&
$\epsilon_{\text{base}} \hspace{.4mm} = 4.281 \cdot 10^{-2}$ & $E_{\max} \hspace{5.1mm} = 4729.19$ \\

	 		& $\gamma_{\text{base}} \hspace{.2mm} = 0$ 						&	$E_{\text{RMSE}} \hspace{3mm} = 39.58$ 
&&
$\gamma_{\text{base}} \hspace{.2mm} = 0$ & $E_{\text{RMSE}} \hspace{3mm} = 1104.12$ \\

			& $\lambda_{\text{base}} = 10^{-5}$ 				&	$E_{\max, \text{rel}} \hspace{2.2mm} = 6.95$ &&
$\lambda_{\text{base}} = 10^{-2}$ & $E_{\max, \text{rel}} \hspace{2.2mm} = 116.02$ \\ 

			& $n_{\text{base}} = 142$ 							&	$E_{\text{RMSE}, \text{rel}} = 9.00 \cdot 
10^{-1}$ &&
$n_{\text{base}} = 63$ & $E_{\text{RMSE}, \text{rel}} = 14.83$ \\ \hline

\multirow{4}{*}{\rotatebox[origin=c]{90}{\textbf{stabilized}}}	& $\epsilon_{\text{stab}} \hspace{.4mm} = \epsilon_{\text{base}}$	&	$E_{\max} 
\hspace{5.1mm} = 344.91$	&&
$\epsilon_{\text{stab}} \hspace{.4mm} = \epsilon_{\text{base}}$ & $E_{\max} \hspace{5.1mm} = 234.40$ \\

			& $\gamma_{\text{stab}} \hspace{.2mm} = 0.5$ 						&	$E_{\text{RMSE}} \hspace{3mm} = 35.69$ 	
&& 
$\gamma_{\text{stab}} \hspace{.2mm} = 0.2$ & $E_{\text{RMSE}} \hspace{3mm} = 30.22$  \\

			& $\lambda_{\text{stab}} = 10^{-5}$ 				&	$E_{\max, \text{rel}} \hspace{2.2mm} = 1.79 \cdot 10^{-1}$ && 
$\lambda_{\text{stab}} = 10^{-15}$ & $E_{\max, \text{rel}} \hspace{2.2mm} = 6.77 \cdot 10^{-1}$  \\

			& $n_{\text{stab}} = 690$  							&	$E_{\text{RMSE}, \text{rel}} = 2.26 \cdot 
10^{-2}$ && 
$n_{\text{stab}} = 358$ & $E_{\text{RMSE}, \text{rel}} = 7.97 \cdot 10^{-2}$  \\
\hline
\end{tabular}
\caption{Overview of the selected hyperparameter and the accuracies of the kernel models.}
\label{tab:num_results}
}
\end{table}

\begin{figure}[ht]
\setlength\fwidth{.395\textwidth}
\centering
\begin{tabular}{cc}
%
%
\begin{tikzpicture}

\begin{axis}[%
width=0.951\fwidth,
height=0.75\fwidth,
at={(0\fwidth,0\fwidth)},
scale only axis,
xmin=0,
xmax=1,
ymin=0,
ymax=1000,
axis background/.style={fill=white},
]
\addplot [color=red, forget plot]
  table[row sep=crcr]{%
0	59.8\\
0.1	171.4\\
0.2	281.2\\
0.3	352\\
0.4	391.4\\
0.5	378.2\\
0.6	363.4\\
0.7	360.4\\
0.8	399\\
0.9	398.8\\
1	410.8\\
};
\addplot [color=black, draw=none, mark=x, mark options={solid, black}, forget plot]
  table[row sep=crcr]{%
0	52\\
0.1	137\\
0.2	295\\
0.3	325\\
0.4	351\\
0.5	323\\
0.6	321\\
0.7	301\\
0.8	353\\
0.9	352\\
1	360\\
};
\addplot [color=black, draw=none, mark=x, mark options={solid, black}, forget plot]
  table[row sep=crcr]{%
0	79\\
0.1	191\\
0.2	235\\
0.3	371\\
0.4	407\\
0.5	398\\
0.6	365\\
0.7	347\\
0.8	379\\
0.9	415\\
1	416\\
};
\addplot [color=black, draw=none, mark=x, mark options={solid, black}, forget plot]
  table[row sep=crcr]{%
0	63\\
0.1	198\\
0.2	334\\
0.3	403\\
0.4	483\\
0.5	460\\
0.6	446\\
0.7	473\\
0.8	494\\
0.9	475\\
1	499\\
};
\addplot [color=black, draw=none, mark=x, mark options={solid, black}, forget plot]
  table[row sep=crcr]{%
0	74\\
0.1	186\\
0.2	296\\
0.3	356\\
0.4	370\\
0.5	395\\
0.6	360\\
0.7	385\\
0.8	417\\
0.9	395\\
1	417\\
};
\addplot [color=black, draw=none, mark=x, mark options={solid, black}, forget plot]
  table[row sep=crcr]{%
0	31\\
0.1	145\\
0.2	246\\
0.3	305\\
0.4	346\\
0.5	315\\
0.6	325\\
0.7	296\\
0.8	352\\
0.9	357\\
1	362\\
};
\end{axis}

\begin{axis}[%
width=1.227\fwidth,
height=0.92\fwidth,
at={(-0.16\fwidth,-0.101\fwidth)},
scale only axis,
xmin=0,
xmax=1,
ymin=0,
ymax=1,
axis line style={draw=none},
ticks=none,
axis x line*=bottom,
axis y line*=left
]
\end{axis}
\end{tikzpicture}%
& \input{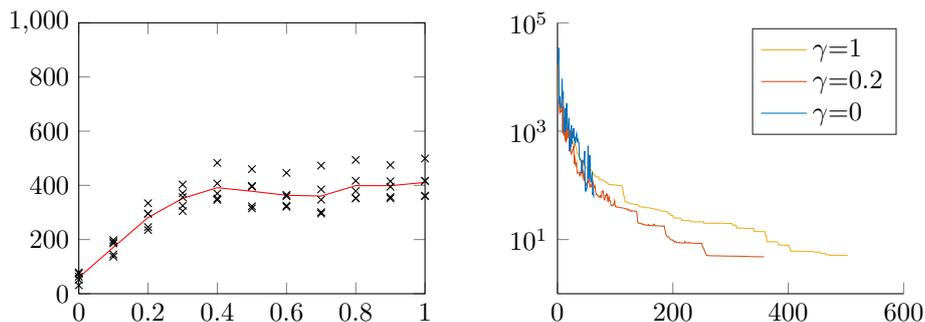}
\end{tabular}

\caption[caption]{Left plot: Number of chosen interpolation points (y-axis) during the 5-fold cross validation procedure for $f/P$-greedy in dependence of the 
restriction parameter $\gamma \in \{0, 0.1, .., 1\}$ (x-axis). The black crosses indicate the five numbers of chosen points during the validation, the red line 
describes the mean value of those. \\ Right plot: $E_{\text{RMSE}}$ error decay (y-axis) during the training of the $f/P$-greedy model depending on the number 
of interpolation points (x-axis) for the unstabilized model ($\gamma=0$), the stabilized model with validated $\gamma$-parameter ($\gamma = 0.2$) and the fully 
stabilized model ($\gamma=1$, i.e. $P$-greedy).}
\label{fig:error_and_chosen_points}
\end{figure}

In the left plot of Figure \ref{fig:error_and_chosen_points} the number of selected points during the cross validation are depicted for the $f/P$-greedy. One 
can see `
that increasing the stabilization parameter $\gamma$ yields more interpolation points. The reason is that the stopping criterion $P_N(x_{N+1}) < \tau_P$ is 
reached later since the selected points are distributed more uniformly as 
quantified in Theorem \ref{th:summary_stabilized} and thus the greedy algorithms run further. We omit plotting results for the $f$-greedy as they do not differ 
considerably. Eventually these further interpolation points yield a better 
interpolation accuracy, which can be seen in Table \ref{tab:num_results}. Especially the maximal relative error $E_{\max, \text{rel}}$ and the relative RMSE 
error $E_{\text{RMSE, rel}}$ are improved. In the right plot of Figure \ref{fig:error_and_chosen_points} the error decay for $f/P$-greedy depending on the 
number of chosen points during the selection (first step of training) is visualized for the $E_{\text{RMSE}}$ error. One can observe that the algorithm stops 
quite early since the stability stopping criterion \eqref{eq:stopping_by_stability} is met. Larger stabilization parameters yield a slower drop, however more 
interpolation points.

\section{Conclusion and Outlook} \label{sec:conclusion}

In this paper a vectorial extension of a recent idea of stabilization of greedy kernel approximation algorithms was introduced and analytical properties were 
stated. A numerical application was addressed using data that emerge in the simulation of the human spine and the stabilization led to significant improvements 
in terms of accuracy and stability due to a better point distribution. \\
A two-step approach was used to combine the stabilization with regularization. In future work we will consider a combined approach of stabilization and 
regularization and use data with more input and output dimensions. Ultimatively we aim at using real patient data and dataset extension approaches by using 
invariances and symmetries or the use of invariant kernels.

\vspace{1cm}
\noindent\textbf{Acknowledgements:} The authors acknowledge the funding of the project by the Deutsche Forschungsgemeinschaft (DFG, German Research Foundation) 
under Germany's Excellence Strategy - EXC 2075 - 390740016.

\bibliography{bibfile_enumath_2019}

\begin{thebibliography}{10}

\bibitem{DeMarchi2005}
S.~De~Marchi, R.~Schaback, and H.~Wendland.
\newblock Near-optimal data-independent point locations for radial basis
  function interpolation.
\newblock {\em Adv. Comput. Math.}, 23(3):317--330, 2005.

\bibitem{DREISCHARF20141757}
M.~Dreischarf, T.~Zander, A.~Shirazi-Adl, C.~Puttlitz, C.~Adam, C.~Chen,
  V.~Goel, A.~Kiapour, Y.~Kim, K.~Labus, J.~Little, W.~Park, Y.~Wang, H.~Wilke,
  A.~Rohlmann, and H.~Schmidt.
\newblock Comparison of eight published static finite element models of the
  intact lumbar spine: Predictive power of models improves when combined
  together.
\newblock {\em Journal of Biomechanics}, 47(8):1757 -- 1766, 2014.

\bibitem{Fasshauer2015}
G.~E. Fasshauer and M.~McCourt.
\newblock {\em Kernel-{B}ased {A}pproximation {M}ethods Using {MATLAB}},
  volume~19 of {\em Interdisciplinary Mathematical Sciences}.
\newblock World Scientific Publishing Co. Pte. Ltd., Hackensack, {NJ}, 2015.

\bibitem{karajan2013a}
N.~Karajan, O.~R{\"o}hrle, W.~Ehlers, and S.~Schmitt.
\newblock Linking continuous and discrete intervertebral disc models through
  homogenisation.
\newblock {\em Biomechanics and Modeling in Mechanobiology}, 12(3):453--66, Jun
  2013.

\bibitem{monteiro2011}
N.~M.~B. Monteiro, M.~P.~T. da~Silva, J.~O. M.~G. Folgado, and J.~P.~L.
  Melancia.
\newblock Structural analysis of the intervertebral discs adjacent to an
  interbody fusion using multibody dynamics and finite element cosimulation.
\newblock {\em Multibody System Dynamics}, 25(2):245--270, Feb 2011.

\bibitem{Mueller2009}
S.~M{\"u}ller.
\newblock {\em Komplexit{\"a}t und Stabilit{\"a}t von kernbasierten
  Rekonstruktionsmethoden (Complexity and Stability of Kernel-based
  Reconstructions)}.
\newblock PhD thesis, Fakult{\"a}t f{\"u}r Mathematik und Informatik,
  Georg-August-Universit{\"a}t G{\"o}ttingen, 2009.

\bibitem{Rupp2015a}
T.~Rupp, W.~Ehlers, N.~Karajan, M.~G{\"u}nther, and S.~Schmitt.
\newblock {A forward dynamics simulation of human lumbar spine flexion
  predicting the load sharing of intervertebral discs, ligaments, and muscles}.
\newblock {\em {Biomechanics and Modeling in Mechanobiology}},
  14(5):1081--1105, 2015.

\bibitem{SH_bookchapter}
G.~Santin and B.~Haasdonk.
\newblock Kernel methods for surrogate modelling.
\newblock Technical Report arXiv:1907.10556, University of Stuttgart, 2019.
\newblock to appear in the MOR Handbook, de Gruyter.

\bibitem{SchWen2000}
R.~Schaback and H.~Wendland.
\newblock Adaptive greedy techniques for approximate solution of large {RBF}
  systems.
\newblock {\em Numer. Algorithms}, 24(3):239--254, 2000.

\bibitem{Wendland2005}
H.~Wendland.
\newblock {\em Scattered {D}ata {A}pproximation}, volume~17 of {\em Cambridge
  Monographs on Applied and Computational Mathematics}.
\newblock Cambridge University Press, Cambridge, 2005.

\bibitem{wenzel2019novel}
T.~{Wenzel}, G.~{Santin}, and B.~{Haasdonk}.
\newblock {A novel class of stabilized greedy kernel approximation algorithms:
  Convergence, stability \&amp; uniform point distribution}.
\newblock {\em arXiv e-prints}, page arXiv:1911.04352, Nov 2019.

\bibitem{Wirtz2013}
D.~Wirtz and B.~Haasdonk.
\newblock A vectorial kernel orthogonal greedy algorithm.
\newblock {\em Dolomites Res. Notes Approx.}, 6:83--100, 2013.

\bibitem{wirtz2015}
D.~Wirtz, N.~Karajan, and B.~Haasdonk.
\newblock Surrogate modeling of multiscale models using kernel methods.
\newblock {\em International Journal for Numerical Methods in Engineering},
  101, 01 2015.

\bibitem{Wittwar2018}
D.~Wittwar, G.~Santin, and B.~Haasdonk.
\newblock Interpolation with uncoupled separable matrix-valued kernels.
\newblock {\em Dolomites Res. Notes Approx.}, 11:23--29, 2018.

\end{thebibliography}
\bibliographystyle{abbrv}

\end{document}